\documentclass[12pt]{article}
\usepackage{natbib,graphicx,setspace,lscape,longtable}
\usepackage{natbib,epsfig,graphicx,float}
\usepackage{amsmath,amsthm,amssymb,color}

\usepackage[nohead]{geometry}
\usepackage[bottom]{footmisc}
\usepackage{indentfirst}
\usepackage{endnotes}
\usepackage{verbatim}
\usepackage{threeparttable}
\usepackage{comment}
\usepackage{algorithm}
\usepackage{titlesec}
\usepackage{multibib}
\usepackage{algorithmic}
\usepackage{multirow}
\usepackage{dsfont}
\usepackage{titling}
\usepackage{afterpage}
\usepackage{hyperref}
\RequirePackage[mathlines, displaymath]{lineno}
\usepackage{multirow}
\usepackage{chngcntr}
\usepackage[titletoc]{appendix}
\bibpunct{(}{)}{;}{a}{,}{,}
\floatstyle{ruled}
\newfloat{algorithm}{tbhp}{loa}
\floatname{algorithm}{Algorithm}

\newcommand{\bA}{\mathrm{\bf A}}

\newcommand{\bZ}{\mathrm{\bf Z}}

\newcommand{\bI}{\mathrm{\bf I}}

\newcommand{\bW}{\mathrm{\bf W}}

\newcommand{\bU}{\mathrm{\bf U}}
\newcommand{\bX}{\mathrm{\bf X}}
\newcommand{\bY}{\mathrm{\bf Y}}

\newcommand{\bzero}{\mathrm{\bf 0}}

\newcommand{\bSigma}{\boldsymbol{\Sigma}}
\newcommand{\bOmega}{\mbox{\boldmath $\Omega$}}

\newcommand{\hSig}{\widehat\Sig}

\newcommand{\Sig}{\mathbf{\Sigma}}

\newcommand{\tr}{\mathrm{tr}}

\newcommand{\var}{\mathrm{var}}

\newcommand{\beq}{\begin{eqnarray*}}
\newcommand{\eeq}{\end{eqnarray*}}

\newcites{New}{References}
\titleformat{\section}{\normalfont\Large\bfseries}{\thesection}{0.5em}{}
\titlespacing*{\section} {0pt}{5pt}{3pt}
\titlespacing*{\subsection} {0pt}{5pt}{2pt}
\numberwithin{equation}{section}
\theoremstyle{plain}

\newtheorem{prop}{Proposition}[section]

\theoremstyle{definition}

\newtheorem{remark}{Remark}[section]
\def\ben{\begin{equation*}}
\def\een{\end{equation*}}
\def\bea{\begin{eqnarray}}
\def\eea{\end{eqnarray}}
\def\bean{\begin{eqnarray*}}
\def\eean{\end{eqnarray*}}
\def\bep{\begin{prop}}
\def\eep{\end{prop}}
\def\bc{\begin{center}}
\def\ec{\end{center}}

\def \tr {\mbox{tr}}
\def \card{\mbox{card}}

\def\lowcomma{_{\textstyle,}}
\def\lowperiod{_{\textstyle.}}

\newtheorem{theorem}{\bf Theorem}
\newtheorem{lemma}{\bf Lemma}

\newtheorem{assump}{\bf Assumption}

\allowdisplaybreaks[4]
\begin{document}
\title{Fisher's combined probability test for high-dimensional covariance matrices\thanks{The authors thank the seminar/conference participants at the George Washington University, University of Southern California, Rutgers University, Yale University, ENAR 2019 Spring Meeting, and 2019 Joint Statistical Meetings for their helpful comments and suggestions. The preliminary result of this paper was included in the National Science Foundation (NSF) grant proposal (DMS–-1811552). Lingzhou Xue's research is supported in part by the NSF grants DMS–-1811552 and DMS-–1953189.}}
\date{First Version: May 2019; This Version: May 2020}
\author{Xiufan Yu, Danning Li, and Lingzhou Xue \\ Pennsylvania State University}
\maketitle{}
\pagestyle{plain}

\begin{abstract}
Testing large covariance matrices is of fundamental importance in statistical analysis with high-dimensional data. In the past decade, three types of test statistics have been studied in the literature: quadratic form statistics, maximum form statistics, and their weighted combination. It is known that quadratic form statistics would suffer from low power against sparse alternatives and maximum form statistics would suffer from low power against dense alternatives. The weighted combination methods were introduced to enhance the power of quadratic form statistics or maximum form statistics when the weights are appropriately chosen.
In this paper, we provide a new perspective to exploit the full potential of quadratic form statistics and maximum form statistics for testing high-dimensional covariance matrices. We propose a scale-invariant power enhancement test based on Fisher's method to combine the $p$-values of quadratic form statistics and maximum form statistics.
After carefully studying the asymptotic joint distribution of quadratic form statistics and maximum form statistics, we prove that the proposed combination method retains the correct asymptotic size and boosts the power against more general alternatives. Moreover, we demonstrate the finite-sample performance  in  simulation studies and a real application.
\end{abstract}


\newpage


\section{Introduction}
Hypothesis testing on large covariance matrices has received considerable attention in the past decade.
The covariance matrices not only have the fundamental importance in multivariate statistics  such as discriminant analysis, principal component analysis, and clustering \citep{anderson2003introduction}, but also play a vital role in various research topics in biological science, finance, operations research including portfolio allocation \citep{goldfarb2003robust}, gene-set testing \citep{chen2010two}, and gene-set clustering \citep{chang2017comparing}.

Let $\bX$ and $\bY$ represent two independent $p$-dimensional random vectors with covariance matrices $\bSigma_1$ and $\bSigma_2$ respectively. We are interested in testing whether these two covariance matrices are equal, that is,
$
H_0:\bSigma_1 = \bSigma_2.
$
This test is well studied in the classical setting where the dimension is fixed and the sample size diverges \citep{anderson2003introduction}. For instance, the likelihood ratio test was shown to enjoy the optimality under mild conditions  \citep*{sugiura1968unbiasedness,perlman1980unbiasedness}.
However, the likelihood function is not well-defined due to the singular sample covariance matrix in the high-dimensional setting where the dimension is no longer fixed but diverges at a possibly faster rate than the sample size.

Over the past decade, statisticians have made a lot of efforts to tackle the challenges in the high-dimensional setting and proposed three different types of statistics for testing large covariance matrices. Firstly, quadratic form statistics were studied to test against the dense alternatives, which can be written in terms of the Frobenius norm of $\bSigma_1-\bSigma_2$ with many small differences between two covariance matrices. When the dimension is on the same order of the sample size, \cite{schott2007test} proposed a test statistic based on the sum of squared differences between two sample covariance matrices, and \cite{srivastava2010testing} used a consistent estimator of $\tr(\bSigma_1^2)/\left[\tr(\bSigma_1)\right]^2 - \tr(\bSigma_2^2)/\left[\tr(\bSigma_2)\right]^2$ to construct a new test statistic.
\cite{li2012two} introduced an unbiased estimator of the Frobenius norm of $\bSigma_1-\bSigma_2$ to allow for the ultra-high dimensionality that the dimension grows much faster than the sample size.
Recently, \cite{he2018asymptotically} proposed the adaptive testing to combine the finite-order U-statistics that includes the variants of quadratic form statistics. Secondly, maximum form statistics were explored to account for the sparse alternatives with only a few large differences between two covariance matrices, which can be written in terms of the entry-wise maximum norm of $\bSigma_1-\bSigma_2$. \cite{cai2013two} studied the maximal standardized differences between two sample covariance matrices to test against the sparse alternative, and \cite{chang2017comparing} proposed a perturbed-based maximum test using a data-driven approach to determine the rejection region.
Thirdly, \cite{li2015joint}, \cite{yang2017weighted} and \cite{li2018applications} used  a  weighted combination of quadratic form statistics and maximum form statistics to test against the dense or sparse alternatives, which shares the similar philosophy with the power enhancement method \citep{fan2015power} for testing cross-sectional dependence.

Similar to these weighted combination tests, we are motivated by combining the strengths of quadratic form statistics and maximum form statistics to boost the power against the dense or sparse alternatives. It is also of great importance to combine the power of these two different statistics in real-world applications such as financial studies and genetic association studies. For instance, the anomalies in financial markets may come from the mispricing of a few assets or a systematic market mispricing \citep{fan2015power}, and the phenotype may be affected by a few causal variants or a large number of mutants \citep{liu2019acat}.

It is worth pointing out that these weighted combination tests critically depend on the proper choice of weights to combine two different types of test statistics. There may exist a non-negligible discrepancy on the different magnitudes between quadratic form statistics and maximum form statistics in practice, which makes the choice of weights a very challenging task. As a promising alternative to \cite{fan2015power}, \cite{li2015joint}, \cite{yang2017weighted} and \cite{li2018applications}, we provide a  new  perspective to exploit the full potential of quadratic form statistics and maximum form statistics for testing high-dimensional covariance matrices.

We propose a scale-invariant power enhancement test based on Fisher's method \citep{Fisher1925} to combine the $p$-values of quadratic form statistics and maximum form statistics.
To study the asymptotic property, we need to solve several non-trivial challenges in the theoretical analysis and then derive the asymptotic joint distribution of quadratic form statistics and maximum form statistics under the null hypothesis. We prove that the asymptotic null distribution of the proposed combination test statistic does not depend on the unknown parameters. More specifically, the proposed statistic follows a chi-squared  distribution with $4$ degrees of freedom asymptotically under the null hypothesis.
We also show the consistent asymptotic power against the union of dense alternatives and sparse alternatives, which is more general than the designated alternative in the weighted combination test. It is worth pointing out that Fisher's method achieves the asymptotic optimality with respect to Bahadur relative efficiency. Moreover, we demonstrate the numerical properties in simulation studies and a real application to gene-set testing \citep{dudoit2008multiple, ritchie2015limma}.
In the real application, the proposed test can detect the important gene-sets more effectively, and our findings are supported by biological evidences.

In recent literature, \cite{liu2019cauchy} proposed the Cauchy combination of  $p$-values for testing high-dimensional mean vectors, and \cite{he2018asymptotically} proved the joint normal limiting distribution of finite-order U-statistics with an identity covariance matrix and used the minimum combination of their $p$-values. The methods and theories of \cite{liu2019cauchy} and \cite{he2018asymptotically} do not apply to the more challenging setting for testing two-sample high-dimensional covariance matrices. Specifically, \cite{li2015joint} and \cite{he2018asymptotically} considered the one-sample test for large covariance matrices that $H_0 : \bSigma = \bI$ under the restricted complete independence assumption among entries of $\bX$, and \cite{li2018applications} studied the one-sample test that $H_0 : \bSigma$ is a banded matrix under the Gaussian assumption.  \cite{li2015joint}, \cite{li2018applications}, and \cite{he2018asymptotically} studied the one-sample covariance test and did not prove the asymptotic independence result for testing two-sample covariance matrices. However, it is significantly more challenging to deal with the complicated dependence in the two-sample tests for large covariance matrices.
To the best of our knowledge, our work presents the first proof of the asymptotic independence result of quadratic form statistics and maximum form statistics for testing two-sample covariance matrices, which provides the essential theoretical guarantee for Fisher's method to combine their $p$-values.


In the theoretical analysis, we use a non-trivial decorrelation technique to address the complex nonlinear dependence in high dimensional covariances. Recently, \cite{shi2019linear} used the decorrelation to study the linear hypothesis testing for high-dimensional generalized linear models. But the nonlinear dependence in the two-sample covariance testing is much more challenging than the linear hypothesis testing. Moreover, we develop a new concentration inequality for two-sample degenerate U-statistics of high-dimensional data, which makes  a  separate contribution to the literature. This result is an extension of the concentration inequality for one-sample degenerate U-statistics \citep{arcones1993limit}.

The rest of this paper is organized as follows. After presenting the preliminaries in Section 2, we introduce the Fisher's method for testing two-sample large covariance matrices in Section 3.  Section 4 studies the asymptotic size and asymptotic power, and Section 5 demonstrates the numerical properties in simulation studies. Section 6 evaluates the proposed test in an empirical study on testing gene-sets. Section 7 includes the concluding remarks. The technical details are presented in the supplementary note.

\section{Preliminaries}

Let $\bX$ and $\bY$ be $p$-dimensional random vectors with covariance matrices $\bSigma_1=\left(\sigma_{ij1}\right)_{p\times p}$ and $\bSigma_2=\left(\sigma_{ij2}\right)_{p\times p}$ respectively. Without loss of generality, we assume both $\bX$ and $\bY$ have zero means.  Let $\left\{ \bX_1,\cdots, \bX_{n_1} \right\}$ be independently and identically distributed (\emph{i.i.d.}) random samples of $\bX$, and $\left\{ \bY_1,\cdots, \bY_{n_2} \right\}$ be \emph{i.i.d.} samples of $\bY$ that are independent of $\left\{ \bX_1,\cdots, \bX_{n_1} \right\}$. The problem of interest is to test  whether two covariance matrices are equal, 
\begin{equation}\label{eq: test}
H_0: \bSigma_1=\bSigma_2. 
\end{equation}

We first revisit the quadratic form statistic \citep{li2012two} to test against the dense alternative and the maximum form statistic \citep{cai2013two} to test against the sparse alternative. The dense alternative can be written in terms of the Frobenius norm of $\bSigma_1-\bSigma_2$ and the sparse alternative can be written using the entry-wise maximum norm of $\bSigma_1-\bSigma_2$.


\cite{li2012two} proposed a quadratic-form test after  reformulating the null hypothesis (\ref{eq: test}) into its equivalent form based on the squared Frobenius norm of $\bSigma_1-\bSigma_2$, that is, $$
H_0: \|\bSigma_1-\bSigma_2\|_F^2= 0. $$ To construct the test statistic, given the simple fact that
$$\|\bSigma_1-\bSigma_2\|_F^2=\tr\{(\bSigma_1-\bSigma_2)^2\} = \tr(\bSigma_1^2) + \tr(\bSigma_2^2) -2\tr(\bSigma_1\bSigma_2),$$
\cite{li2012two} proposed a test statistic $T_{n_1,n_2}$ in the form of linear combination of unbiased estimators for each term, specifically,
\begin{equation}
T_{n_1,n_2}=A_{n_1}+B_{n_2}-2C_{n_1,n_2},
\end{equation}
where $A_{n_1}$, $B_{n_2}$ and $C_{n_1,n_2}$ are the unbiased estimators under $H_0$ for $\tr(\bSigma_1^2)$, $\tr(\bSigma_2^2)$ and $\tr(\bSigma_1\bSigma_2)$ respectively. Then, the expected value of $T_{n_1,n_2}$ is zero under the null hypothesis.
For details about $A_{n_1}$, $B_{n_2}$ and $C_{n_1,n_2}$, please refer to  Section 2 of \cite{li2012two}. \cite{li2012two} proved that the asymptotic distribution of $T_{n_1,n_2}$ is a normal distribution. Let $z_\alpha$ be the upper $\alpha$ quantile of the standard normal distribution, and $\widehat\sigma_{0,n_1,n_2}$ is a consistent estimator of  the leading term $\sigma_{0,n_1,n_2}$ in the standard deviation of $T_{n_1,n_2}$ under $H_0$. Hence, \cite{li2012two} rejects the null hypothesis at the significance level $\alpha$ if
\begin{equation}\label{test: LC}
  T_{n_1,n_2}\geq \widehat\sigma_{0,n_1,n_2}z_\alpha.
\end{equation}

As an alternative to the quadratic form statistic \citep{li2012two},
\cite{cai2013two}
studied the null hypothesis (\ref{eq: test}) in terms of the maximal absolute difference of two covariance matrices, i.e.,
$$
H_0: 
\max_{1\leq i\leq j\leq p} |\sigma_{ij1}-\sigma_{ij2}| = 0. $$
\cite{cai2013two}  proposed a maximum test statistic $M_{n_1,n_2}$ based on the maximum of standardized differences between  $\widehat\sigma_{ij1}$'s and $\widehat\sigma_{ij2}$'s. The maximum form statistic is written as
\begin{equation}
M_{n_1,n_2}=\max_{1\leq i \leq j\leq p} \frac{\left( \widehat{\sigma}_{ij1}-\widehat{\sigma}_{ij2} \right)^2}{\widehat{\theta}_{ij1}/n_1+\widehat{\theta}_{ij2}/n_2}\lowcomma
\end{equation}
where the denominator $\widehat\theta_{ij1}/n_1+\widehat\theta_{ij1}/n_2$ estimates the variance of $\widehat\sigma_{ij1}-\widehat\sigma_{ij2}$ to account for the heteroscedasticity of $\widehat\sigma_{ij1}$'s and $\widehat\sigma_{ij2}$'s among different entries.
\cite{cai2013two} proved that the asymptotic null distribution of $M_{n_1,n_2}$ is a Type \uppercase\expandafter{\romannumeral1} extreme value distribution (also known as the Gumbel distribution). Thus, \cite{cai2013two} rejects the null hypothesis at a significance level $\alpha$ if
\begin{equation}\label{test: CLX}
M_{n_1,n_2} \geq q_{\alpha}+4\log p- \log\log p\lowcomma
\end{equation}
where $q_\alpha$ is the upper $\alpha$ quantile of the Gumbel distribution.

\section{Fisher's Combined Probability Test}

\cite{li2012two} and \cite{cai2013two} have their respective power for testing high-dimensional covariance matrices. The quadratic form statistic $T_{n_1,n_2}$ is powerful against the dense alternative, where the difference between $\bSigma_1$ and $\bSigma_2$ under the squared Frobenius norm is no smaller than the order of $\tr(\bSigma_1^2)/n_1+\tr(\bSigma_2^2)/n_2$.
The maximum form statistic  $M_{n_1,n_2}$ is powerful against the sparse alternative, where at least one entry of $\bSigma_1-\bSigma_2$ has the magnitude larger than the order of $\sqrt{\log p/n}$. However, $T_{n_1,n_2}$ performs poorly against the sparse alternative and $M_{n_1,n_2}$ performs poorly against the dense alternative. More details will be presented in Subsection \ref{subsec: size-power} and Section \ref{sec:simulation}.

\cite{fan2015power}, \cite{li2015joint}, \cite{yang2017weighted} and \cite{li2018applications} studied the weighted combination $J = J_0+J_1$ to achieve the power enhancement, where $J_0$ is built on the extreme value form statistic and $J_1$ is constructed from the asymptotically pivotal statistic. It is worth pointing out that, with the proper weighted combination, $J$ enjoys the so-called \textsl{power enhancement principles} \citep{fan2015power}: (i) $J$ is at least as powerful as $J_1$, (ii) the size distortion due to the addition of $J_0$ is asymptotically negligible, and (iii) power is improved under the designated alternatives. For testing large covariance matrices, \cite{yang2017weighted} proposed $J_1 = (1-(s_p+\xi_1)^{-1})M_n$ and $J_0 = n^{\frac{1}{s_p+\xi_1} + \frac{1}{\xi_2}}\cdot \max_{1\leq i,j\leq p} (\widehat\sigma_{ij1}-\widehat\sigma_{ij2})^2 $, where $M_n$ is a macro-statistic  which performs well against the dense alternative, and $s_p$ is the number of distinct entries in two covariance matrices. Note that the quantities $\xi_1$ and $\xi_2$ are carefully chosen such that $J_0 \rightarrow 0$ under $H_0$.


As a promising alternative, we propose a scale-invariant combination procedure based on Fisher's method \citep{Fisher1925} to combine both strengths of $T_{n_1,n_2}$ and $M_{n_1,n_2}$. Let $\Phi(\cdot)$ be the cumulative distribution function of $N(0,1)$ and $G(x)=\exp\left(-\frac{1}{\sqrt{8\pi}}\exp\left(-\frac{x}{2}\right)\right)$ be the cumulative distribution function of the Gumbel distribution. More specifically, we combine the $p$-values of $T_{n_1,n_2}$ and $M_{n_1,n_2}$ after the negative natural logarithm transformation, that is,
\begin{equation}
F_{n_1,n_2} = -2 \log p_T - 2 \log p_M,
\end{equation}
where $$p_T = 1-\Phi\left(T_{n_1,n_2}/\widehat\sigma_{0,n_1,n_2}\right)$$ and $$p_M = 1-G(M_{n_1,n_2}-4\log p+\log\log p)$$ are the $p$-values associated with the test statistics $T_{n_1,n_2}$ and $M_{n_1,n_2}$, respectively.


Let $c_\alpha$ denote the upper $\alpha$ quantile of a  chi-squared  distribution  with  4  degrees  of  freedom (i.e., $\chi_4^2$). We reject the null hypothesis at the significance level $\alpha$ if
\begin{equation}
F_{n_1,n_2}\geq c_\alpha.
\end{equation}

Unlike the weighted statistic $J = J_0+J_1$, $F_{n_1,n_2}$ does not need to estimate $s_p$ or choose $\xi_1$ and $\xi_2$ to construct the proper weights, which may be non-trivial to deal with in practice.
The inappropriate choice of $s_p$, $\xi_1$ and $\xi_2$ may lead to the  size distortion or loss of power. In contrast, $F_{n_1,n_2}$ is scale-invariant as the $p$-values always take values between 0 and 1, and the asymptotic null distribution of $F_{n_1,n_2}$ (i.e., $\chi_4^2$) does not depend on any hyper-parameters. As we will show in Section \ref{subsec: size-power}, $F_{n_1,n_2}$ achieves the desired nominal significance level asymptotically while boosting the power against either sparse or dense alternatives. Moreover, Fisher's method achieves the asymptotic optimality with respect to Bahadur relative efficiency \citep{littell1971asymptotic,littell1973asymptotic}.

\begin{remark}
The idea of combining $p$-values has been widely used as an important technique for data fusion or meta analysis \citep{hedges2014statistical}. Recently, the Cauchy combination of  $p$-values was used for testing high-dimensional mean vectors in \citep{liu2019cauchy}, and the minimum combination of $p$-values from the finite-order U-statistics was used for testing two-sample high-dimensional covariance matrices in \citep{he2018asymptotically}. However, neither \cite{liu2019cauchy} nor \cite{he2018asymptotically} studied the combination of $p$-values of $T_{n_1,n_2}$ and $M_{n_1,n_2}$, 
and it is fundamentally challenging to study the asymptotic joint distribution of $T_{n_1,n_2}$ and $M_{n_1,n_2}$. We will solve this open problem in Subsection 4.2.
\end{remark}

\section{Asymptotic Properties}\label{sec: asymptotic}
This section presents the asymptotic properties of our proposed Fisher's combined probability test $F_{n_1,n_2}$. Section \ref{subsec: assumptions} presents the assumptions. Section \ref{subsec: asymp-independence} studies the joint limiting distribution of two test statistics $M_{n_1,n_2}$ and $T_{n_1,n_2}$ under the null hypothesis. Section
\ref{subsec: size-power} proves the correct asymptotic size and consistent asymptotic power of our proposed method.
\subsection{Assumptions}\label{subsec: assumptions}

We define some useful notations. For any matrix $\bA$, let $\lambda_i(\bA)$ be the $i$-th largest eigenvalue of $\bA$. For any set $\mathcal{A}$, $\card(\mathcal{A})$ represents the cardinality of $\mathcal{A}$. For  $0<r<1$, let
\begin{equation*}
\mathcal{V}_i(r) = \left\{1\leq j\leq p: \frac{|\sigma_{ij1}|}{\sqrt{\sigma_{ii1}\sigma_{jj2}}}\geq r \text{ or } \frac{|\sigma_{ij2}|}{\sqrt{\sigma_{ii2}\sigma_{jj2}}} \geq r \right\}
\end{equation*}
be the set of indices $j$ such that $X_j$ (or $Y_j$) is highly correlated (whose correlation $>r$) with $X_i$ (or $Y_i$) for a given $i\in\{1,\dots, p\}$. And for any $\alpha>0$, let
\begin{equation*}
    s_i(\alpha) = \card(\mathcal{V}_i(\left(\log p\right)^{-1- \alpha})),\ i=1,\cdots, p
\end{equation*}
denote the number of indices $j$ in the set $\mathcal{V}_i(\left(\log p\right)^{-1- \alpha})$. Moreover, define
\begin{equation*}
    \mathcal{W}(r) = \left\{1\leq i\leq p: \mathcal{V}_i(r) \neq \varnothing \right\}
\end{equation*}
such that, $\forall i\in\mathcal{W}(r)$, $X_i$ (or $Y_i$) is highly correlated with some other variable of $\bX$ (or $\bY$).

Throughout the rest of this section, we assume that $\bX$ and $\bY$ are both Gaussian random vectors. The Gaussian assumption facilitates the use of a new decorrelation technique to address the complex nonlinear dependence in high dimensional covariances in the theoretical analysis of the proposed scale-invariant combination test.

\begin{remark}
\cite{li2015joint}, \cite{li2018applications} and \cite{he2018asymptotically} studied the asymptotic joint distribution of the maximum test statistic and the quadratic test statistic for one-sample covariance test under the Gaussian assumption or restricted
complete independence assumption. Please see the first paragraph of Section 2 in \cite{li2015joint}, the first paragraph of Section 2 in \cite{li2018applications}, and Condition 2.3 in \cite{he2018asymptotically} for more details. However, the nonlinear dependence in two-sample covariance test is fundamentally more challenging than the dependence in the  one-sample covariance test.
\end{remark}

\begin{assump}\label{assum: A1A2-in-Chen}
As $\min\{n_1,n_2\}\rightarrow\infty$ and $p \rightarrow\infty$,
\begin{itemize}
  \item[(i)] $n_1/\left(n_1+n_2\right)\rightarrow \gamma$, for some constant $\gamma\in (0,1)$.
  \item[(ii)] $\sum_{i=1}^{q}\lambda^2_i(\bSigma_j)/\sum_{i=1}^{p}\lambda^2_i(\bSigma_j)\to 0$ for any integer $q = O(\log p)$ and $j=1, 2$.
\end{itemize}

\end{assump}

\begin{remark}
Assumption \ref{assum: A1A2-in-Chen} is analogous to (A1) and (A2) in \cite{li2012two}, where the first condition is standard for two-sample asymptotic analysis, and the second one describes the extent of high dimensionality and the dependence which can be accommodated by the proposed tests. Sharing the spirit, Assumption \ref{assum: A1A2-in-Chen} does not impose explicit requirements on relationships between $p$ and $n_1, n_2$, but rather
requires a mild condition (ii) regarding the covariances, which can be satisfied if eigenvalues of two covariance matrices are bounded.
\end{remark}

\begin{assump}\label{assum: C1-in-Cai}
There exists a subset $\Upsilon\subset\left\{1,2, \cdots,p\right\}$ with $\card\left(\Upsilon\right)=o(p)$ and some constant $\alpha_0>0$, such that for all $\kappa>0$, $\operatornamewithlimits{\max}\limits_{1\leq i\leq p, i\not\in \Upsilon} s_i(\alpha_0)=o(p^\kappa)$. 
In addition, there exists a constant $0<r_0<1$, such that $\card(\mathcal{W}(r_0)) = o(p)$.

\end{assump}

\begin{remark}
Assumption \ref{assum: C1-in-Cai} was introduced by \cite{cai2013two} such that $\max_{1\leq i\leq p, i\not\in \Upsilon}$ $s_i(\alpha_0)$ and $\mathcal{W}(r_0)$ are moderate for $\alpha_0>0$ and $0<r_0<1$. It is satisfied if the eigenvalues of covariance matrices are bounded from above and correlations are bounded away from $\pm 1$.
\end{remark}

\subsection{Asymptotic Joint Distribution}\label{subsec: asymp-independence}
Now, we present the joint limiting law for $M_{n_1,n_2}$ and $T_{n_1,n_2}$ under the null hypothesis.

\begin{theorem}\label{thm: asymp-indep}
Suppose Assumptions \ref{assum: A1A2-in-Chen} and  \ref{assum: C1-in-Cai} hold, and $\log p = o(n^{\frac{1}{5}})$ for $n = n_1 + n_2$, then under the null hypothesis $H_0$, for any $x,t\in\mathbb{R}$, we have
\begin{equation}\label{eq: asympindep}
P\left(\frac{T_{n_1,n_2}}{\widehat\sigma_{0,n_1,n_2}}\leq t,\ M_{n_1,n_2}-4\log p+\log\log p \leq x\right) \rightarrow \Phi(t)\cdot G(x)
\end{equation}
as $n_1,n_2,p\rightarrow\infty$, where $G(x)=\exp\left(-\frac{1}{\sqrt{8\pi}}\exp\left(-\frac{x}{2}\right)\right)$ is the cdf of Gumbel distribution, and $\Phi(t)$ is the cdf of standard normal distribution.
\end{theorem}

\begin{remark}
Together with Theorems 1 and 2 from \cite{li2012two} and Theorem 1 from \cite{cai2013two}, Theorem \ref{thm: asymp-indep} implies that $M_{n_1,n_2}$ and $T_{n_1,n_2}$ are asymptotically independent.
\end{remark}

In the sequel, we provide a high-level intuition to prove the asymptotic independence result (\ref{eq: asympindep}). First of all, it is worth mentioning that under Assumption \ref{assum: A1A2-in-Chen},  all the third-moment and fourth-moment terms in $A_{n_1}$, $B_{n_2}$ and $C_{n_1,n_2}$ are of small order than the leading second-moment terms, which may be neglected when deriving the asymptotic normality.
Hence in theoretical analysis, we may consider the simplified  statistic of $T_{n_1,n_2}$ defined by
\begin{equation}\label{eq: Tp}
\widetilde{T}_{n_1,n_2}=\frac{1}{n_1(n_1-1)}\sum_{u\neq v}\left(\bX_u'\bX_v\right)^2+\frac{1}{n_2(n_2-1)}\sum_{u\neq v}\left(\bY_u'\bY_v\right)^2-\frac{2}{n_1n_2}\sum_u\sum_v\left(\bX_u'\bY_v\right)^2\lowperiod
\end{equation}
As pointed out by \cite{li2012two}, $\widetilde{T}_{n_1,n_2}$ and $T_{n_1,n_2}$ shares the same asymptotic behavior.

Compared with the simple one-sample covariance test in \cite{li2015joint}, \cite{li2018applications}, and \cite{he2018asymptotically}, it is significantly more difficult to analyze the asymptotic joint distribution given the complicated dependence in the two-sample tests for large covariance matrices. To address this challenge, we use a decorrelation technique to address the complex nonlinear dependence in high dimensional covariances. Specifically, we introduce a decorrelated statistic $T_{n_1,n_2}^{*}$. Under  $H_0: \bSigma_1 = \bSigma_{2} = \bSigma$, we may partition $\bX$ and  $\bY$ as follows:
\begin{equation*}
\bX_{p\times 1}=\begin{pmatrix} \bX^{(1)} \\ \bX^{(2)} \end{pmatrix} \text{ and } \bY_{p\times 1}=\begin{pmatrix} \bY^{(1)} \\ \bY^{(2)} \end{pmatrix} 
\sim N_{p}\left( \begin{pmatrix}  \bzero_{p-q} \\  \bzero_{q} \end{pmatrix}\lowcomma \bSigma=\begin{pmatrix} \bSigma_{11} & \bSigma_{12} \\ \bSigma_{21} & \bSigma_{22} \end{pmatrix} \right).
\end{equation*}
where $\bX^{(1)}, \bY^{(1)}\in\mathbb{R}^{p-q},\ \bX^{(2)}, \bY^{(2)}\in\mathbb{R}^{q}$ for integer $q$ satisfying $q=O(\log p)$.
Let $\bZ_1=\bX^{(1)}-\bSigma_{12}\bSigma_{22}^{-1}\bX^{(2)}$, $\bZ_2=\bX^{(2)}$, $\bW_1=\bY^{(1)}-\bSigma_{12}\bSigma_{22}^{-1}\bY^{(2)}$, $\bW_2=\bY^{(2)}$. It's easy to see that $\bZ_1$ is independent of $\bZ_2$, and the same results hold for $\bW_1$ and $\bW_2$. Back to the sample level, we have that $\{\bZ_{1u}\}_{u=1}^{n_1}$ and $\{\bW_{1v}\}_{v=1}^{n_2}$ i.i.d. follow $N_{p-q}(\bzero, \bSigma_{11}-\bSigma_{12}\bSigma_{22}^{-1}\bSigma_{21})$. Following the pattern of $\widetilde{T}_{n_1,n_2}$ in (\ref{eq: Tp}), we define
\begin{equation}\label{eq: Tstar}
T_{n_1,n_2}^*=\frac{1}{n_1(n_1-1)}\sum_{u\neq v}\left(\bZ_{1u}'\bZ_{1v}\right)^2+\frac{1}{n_2(n_2-1)}\sum_{u\neq v}\left(\bW_{1u}'\bW_{1v}\right)^2-\frac{2}{n_1n_2}\sum_u\sum_v\left(\bZ_{1u}'\bW_{1v}\right)^2\lowperiod
\end{equation}
$\{\bZ_{1u}\}_{u=1}^{n_1}$ and $\{\bW_{1v}\}_{v=1}^{n_2}$ are regarded as a decorrelated version of $\{\bX_{u}\}_{u=1}^{n_1}$ and $\{\bY_{v}\}_{v=1}^{n_2}$, respectively. $T_{n_1,n_2}^*$ is regarded as the $\widetilde{T}_{n_1,n_2}$ statistic derived from the decorrelated samples. The above deccorelation shares a similar philosophy with \citep{shi2019linear}. We should point out that \cite{shi2019linear} used the decorrelation to study the linear hypothesis testing for high-dimensional generalized linear models, but the nonlinear dependence in the two-sample covariance testing is much more challenging than the linear hypothesis testing.

In what follows, we study the joint distribution of $M_{n_1,n_2}$ and $\widetilde{T}_{n_1, n_2}$. Let $A$ denote the event associated with the maximum statistic $M_{n_1, n_2}$, and let $B$ be the event corresponding to the quadratic statistic $\widetilde{T}_{n_1,n_2}$. We use the simple but very helpful fact that $A = \cup_i A_i$. Then, we may rewrite the joint probability $P\left(A\cap B\right)$ into the probability for a union of events, that is, $P\left(A\cap B\right) = P\left( (\cup_i A_i)\cap B\right)$. In what follows, we give the proof sketch to derive the upper bound $P(A\cap B) - P(A)P(B)\le o(1)$. We begin with a union bound to obtain that $P\left(\cup_i(A_i\cap B)\right)\leq \sum_{i} P(A_i\cap B)$. In order to deal with the joint probability of $A_{i} \cap B$, we further decompose the quadratic statistic into two parts: $T_{n_1,n_2}^*$ is independent of $A_i$, and the remaining term $ \widetilde{T}_{n_1,n_2}- T_{n_1,n_2}^*$ is associated with $A_i$. Consequently, $B$ can be written as $B = B_i^c \cup B_i$, in which $B_i^c$ represents to the event corresponding to $T_{n_1,n_2}^{*}$. Therefore,  $\sum_i P(A_i\cap B) \leq \sum_{i} P(A_i\cap B_i^c) + \sum_i P(A_i\cap B_i)\leq \sum_{i} P(A_i)P(B_i^c) +\sum_{i} P(B_i)$.
Lemma \ref{lem: T-star-expdecay} suggests $T_{n_1,n_2}^*$ is sufficiently close to $\widetilde{T}_{n_1,n_2}$ so that we have $P(B_i^c)\approx P(B)$, $\sum_{i} P(A_i)\to P(A)$ and $\sum_{i} P(B_i) = o(1)$.
The lower bound $o(1) \le P(A\cap B) - P(A)P(B)$ can be similarly derived from the Bonferroni inequality. Therefore, we can prove the asymptotic independence given that $|P(A\cap B) - P(A)P(B)| = o(1)$.


In the following, we present three useful lemmas to prove (\ref{eq: asympindep}) in Theorem 1.

\begin{lemma}[Asymptotic Normality]\label{lem: T-star-asympnormality} Under Assumption \ref{assum: A1A2-in-Chen},  as $n_1,n_2,p\rightarrow\infty$,
\begin{equation}\label{eq: T-star}
\frac{T_{n_1,n_2}^*}{2\left(n_1^{-1}+n_2^{-1}\right)\tr\left(\bSigma^2\right)} \overset{d}{\rightarrow} N(0,1).
\end{equation}
\end{lemma}

\begin{lemma}[Exponential Decay]\label{lem: T-star-expdecay} Under Assumption \ref{assum: A1A2-in-Chen}, for any $\epsilon>0$, there exists positive constants $C, c$ that do not depend on $p$, $n_1$, $n_2$, such that
\begin{equation}\label{eq: exp-decay}
P\left( \frac{\left|\widetilde{T}_{n_1,n_2}-T_{n_1,n_2}^*\right|}{2\left(n_1^{-1}+n_2^{-1}\right)\tr\left(\bSigma^2\right)} \geq \epsilon \right)    \leq C \exp\{-c \epsilon n^{\beta}\},
\end{equation}
with $1/5<\beta<1/3.$
\end{lemma}

\begin{remark}
Lemma \ref{lem: T-star-expdecay} presents a new concentration inequality for two-sample degenerate U-statistics. It extends the well-known concentration inequality for one-sample degenerate U-statistics \citep{arcones1993limit} and makes  a  separate contribution to the literature.
\end{remark}

As a final step, Lemma \ref{lem: asymp-indep} derives the joint limiting distribution of the test statistic $M_{n_1, n_2}$ and the simplified statistic $\widetilde{T}_{n_1,n_2}$, which directly implies Theorem \ref{thm: asymp-indep}.

\begin{lemma}\label{lem: asymp-indep}
Under the same assumptions as in Theorem \ref{thm: asymp-indep},
\begin{equation}\label{eq: asympindep2}
P\left(\frac{\widetilde{T}_{n_1,n_2}}{\widehat\sigma_{0,n_1,n_2}}\leq t,\ M_{n_1,n_2}-4\log p+\log\log p \leq x\right) {\rightarrow} \Phi(t)\cdot G(x)
\end{equation}
for any $x,t\in\mathbb{R}$, as $n_1,n_2,p\rightarrow\infty$.
\end{lemma}

 Lemma \ref{lem: T-star-asympnormality} shows that such decorrelation procedure does not affect the asymptotic behavior of the quadratic test statistic. Lemma \ref{lem: T-star-expdecay} depicts the tail behavior of the difference between $\widetilde{T}_{n_1,n_2}$ and ${T}_{n_1,n_2}^*$ with explicit decaying rate.  Lemma \ref{lem: T-star-asympnormality} and Lemma \ref{lem: T-star-expdecay} lay the foundation of replacing $\widetilde{T}_{n_1,n_2}$ with  ${T}_{n_1,n_2}^*$ in the theoretical analysis.

\subsection{Asymptotic Size and Power}\label{subsec: size-power}

Given the explicit joint distribution of $M_{n_1, n_2}$ and $T_{n_1,n_2}$, we proceed to present the asymptotic properties of our proposed Fisher's test $F_{n_1,n_2}$. Recall that $c_\alpha$ is the upper $\alpha$-quantile of $\chi_4^2$ distribution and $F_{n_1, n_2} = -2\log (p_M) - 2\log(p_T)$ rejects $H_0$ if $F_{n_1, n_2}$ is as extreme as $c_\alpha$. On top of the asymptotic independence established in Section \ref{subsec: asymp-independence} and by simple probability transformation, it's easy to obtain the null distribution of $F_{n_1, n_2}$, and therefore, the asymptotic size of the test. The results are formally presented in Theorem \ref{thm: size}.

\begin{theorem}[Asymptotic Size]\label{thm: size}
Under the same assumptions as in Theorem \ref{thm: asymp-indep}, the Fisher's test achieves accurate asymptotic size, that is, under the null hypothesis,
$$P\left(F_{n_1,n_2}\geq c_\alpha\right) \rightarrow \alpha\quad \text{as } n_1,n_2,p\rightarrow \infty.$$
\end{theorem}

\begin{remark}
Besides Fisher's method, the asymptotic independence result makes it feasible to combine $p$-values using other approaches such as Tippett's method \citep{tippett1931methods}, Stouffer's method \citep{stouffer1949american}, and Cauchy combination \citep{liu2019cauchy}.
\end{remark}

\cite{li2012two} and \cite{cai2013two} provided power analysis of tests $T_{n_1,n_2}$ and $M_{n_1,n_2}$ over the dense alternative $\mathcal{G}_d$ and the sparse alternative $\mathcal{G}_s$ respectively.
\begin{align}
\mathcal{G}_d & = \left\{ (\bSigma_1,\bSigma_2): \bSigma_1 >0, \bSigma_2>0, \frac{1}{n_1}\tr(\bSigma_1^2) + \frac{1}{n_2}\tr(\bSigma_2^2) = o\left(\tr\{(\bSigma_1-\bSigma_2)^2\}\right)\right\} \label{eq: G1}\lowcomma \\
\mathcal{G}_s & = \left\{(\bSigma_1,\bSigma_2): \bSigma_1 >0, \bSigma_2>0, \max_{1\leq i\leq j\leq p} \frac{|\sigma_{ij1}-\sigma_{ij2}|}{\sqrt{\theta_{ij1}/n_1+\theta_{ij2}/n_2}} \geq 4\sqrt{\log p} \right\}\lowperiod \label{eq: G2}
\end{align}

Taking advantage of the combination, we shall show that our proposed combined test $F_{n_1,n_2}$ makes the most of merits from the two tests and successfully boost the power against either dense or sparse alternatives.

\begin{theorem}[Asymptotic  Power]\label{thm: power}
Under the same assumptions as in Theorem \ref{thm: asymp-indep}, the Fisher's test achieves consistent asymptotic power, that is, under the alternative hypothesis, $$\inf_{(\bSigma_1,\bSigma_2)\in \mathcal{G}_d \cup \mathcal{G}_s} P\left(F_{n_1,n_2}\geq c_\alpha\right) \rightarrow 1 \quad \text{as } n_1,n_2,p\rightarrow \infty.$$

\end{theorem}

\begin{remark}
\noindent {(Bahadur Efficiency)} As discusses in \cite{littell1971asymptotic,littell1973asymptotic}, among all approaches of combining independent tests, Fisher's method delivers the largest exact Bahadur slope, indicating the fastest decay rate of the p-values. Therefore, 
Fisher's test is asymptotically optimal in terms of Bahadur relative efficiency.
\end{remark}

\section{Simulation Studies}\label{sec:simulation}

This section examines the finite-sample performance of our Fisher's combined probability test, compared to the tests proposed by \cite{cai2013two} (refer as the CLX test in the following context) and \cite{li2012two} (refer as the LC test). We generate $\{\bX_1,\cdots,\bX_{n_1}\}$ \emph{i.i.d.} from $N_p\left(\bzero,\bSigma_1\right)$ and $\{\bY_1,\cdots,\bY_{n_2}\}$ \emph{i.i.d.} from $N_{p}\left(\bzero, \bSigma_2\right)$. 
The sample sizes are taken to be $n_1=n_2=N$ with $N=100$ and $200$, while the dimension $p$ varies over the values 100, 200, 500, 800 and 1000. For each simulation setting, the average number of rejections are reported based on 1000 replications. The significance level is set to be $0.05$ for all the tests.

Under the null hypothesis $H_0$, we set $\bSigma_1=\bSigma_2=\bSigma^{*(i)}, i=1,\cdots,5$, and consider the following five models to evaluate the testing size.
\begin{itemize}
	\item[(i)] $\bSigma^{*(1)}=\bI_p$.
	\item[(ii)] $\bSigma^{*(2)}=(\bOmega^{*(2)})^{-1}$, where $\omega_{ij}^{*(2)}=0.5^{|i-j|}$.
	\item[(iii)] $\bSigma^{*(3)}$ is a block diagnoal matrix given by each block being $0.5\bI_5+0.5\mathds{1}_5\mathds{1}'_5$.
	\item[(iv)] $\bSigma^{*(4)}=\{\sigma_{ij}^{*(4)}\}_{p\times p}$, $\sigma_{ij}^{*(4)}=(-1)^{i+j}0.4^{|i-j|^{1/10}}$.
	\item[(v)] $\bSigma^{*(5)}=(\bSigma^{(5)}+\delta\bI)/(1+\delta)$, where $\sigma_{ii}^{(5)}=1$, $\sigma_{ij}^{(5)}=0.5*Bernoulli(1,0.05)$ for $i<j$ and $\sigma_{ij}^{(5)}=\sigma_{ji}^{(5)}$, $\delta=|\lambda_{\min}(\bSigma^{(5)})|+0.05$.
\end{itemize}

Model (i) is the most commonly used multivariate standard normal distribution. Model (ii) and Model (iii) are the cases when the true covariance matrices have certain banded-type and block-type sparsity. Model (iv) was first proposed by \cite{srivastava2010testing} and further studied in \cite{cai2013two}. Model (v) is also a sparse matrix yet without any specific sparsity pattern.

To evaluate the power of the tests, we consider the scenarios when the differences of the two covariance matrices satisfy certain structure. There are two types of alternatives we desire to look into: the sparse alternative $H_s$ and the dense alternative $H_d$.

Generally speaking, the sparse alternative shares commonality among different models. Let $\bU$ denote the difference between $\bSigma_2$ and $\bSigma_1$, i.e. $\bU=\bSigma_2-\bSigma_1$. Inspired by \cite{cai2013two}, we consider the situation when $\bU$ is a symmetric sparse matrix with eight random nonzero entries. The locations of four nonzero entries are randomly selected from the upper triangle of $\bU$, each with a magnitude of Unif(0,4)$\times\max_{1\leq j \leq p} \sigma_{jj}^{*}$. The other four are determined by symmetry. Then we generate samples from  these covariance pairs $\left(\bSigma_{1}^{(i)}, \bSigma_{2}^{(i)}\right)$, $i=1,\cdots,5$, in order to evaluate the power of the tests against sparse alternative, where $\bSigma_{1}^{(i)}=\bSigma^{*(i)}+\delta \bI$ and $\bSigma_{2}^{(i)}=\bSigma^{*(i)}+\delta \bI+ \bU$, with $\delta=|\min\{\lambda_{\min}(\bSigma^{*(i)}+\bU),\lambda_{\min}(\bSigma^{*(i)})\}|+0.05$.

In terms of the dense alternative setting, since the five  models differ a lot from each other, we shall discuss  their corresponding alternative settings separately afterwards. To begin with, we shall take a look at the simplest case in Model (i). We consider its dense alternative to be the AR(1) model with parameter $\rho=0.2$ and $0.3$, denoted by $\bSigma_{\rho}^{AR}$. In another word, we generate the copies of $\bX$ from the $p$-dimensional standard normal while copies of $\bY$ from $N_p\left(\bzero,\bSigma_{\rho}^{AR}\right)$. We follow the same alternative hypothesis as in \cite{srivastava2010testing} for Model (iv), which is $\sigma_{ij}^{(4)}=(-1)^{i+j}0.6^{|i-j|^{1/10}}$, whereas we use the identity matrix $\bI_p$ for Models (ii), (iii) and (v).

\begin{table}[H]
\small
\centering
\caption{Comparison of Empirical Size and Power (\%) for Model (i)}\label{tab: standardnormal}
\vspace{1ex}
\begin{tabular}{ccrrrrr|rrrrrrrrrrr}
\hline
 n&p& 100 & 200 & 500 & 800 & 1000 & 100 & 200 & 500 & 800 & 1000 \\
   \hline
 &&\multicolumn{5}{c|}{Size} &\multicolumn{5}{c}{Power under sparse alternative}\\
\multirow{3}{*}{100}
& Proposed & 5.6 & 5.0 & 5.0 & 5.2 & 5.6 & 98.0 & 96.6 & 87.3 & 83.9 & 80.2\\
& CLX      & 4.3 & 5.2 & 4.5 & 4.4 & 4.5 & 98.5 & 98.3 & 91.1 & 89.8 & 85.8\\
& LC       & 4.8 & 5.0 & 5.1 & 4.5 & 4.2 & 20.6 & 11.2 &  5.9 &  5.7 &  5.0\\
\hline
\multirow{3}{*}{200}
& Proposed & 4.6 & 4.7 & 4.8 & 4.9 & 4.3 & 100.0 & 100.0 & 100.0 & 100.0 & 100.0\\
& CLX      & 3.6 & 4.2 & 4.5 & 5.5 & 5.0 & 100.0 & 100.0 & 100.0 & 100.0 & 100.0 \\
& LC       & 5.4 & 3.2 & 4.6 & 4.8 & 5.3 &  50.5 &  22.2 &   8.0 &   7.6 &   7.3 \\
\hline
 &&\multicolumn{10}{c}{Power under dense alternative}\\
  &&  \multicolumn{5}{c}{$\rho=0.2$} & \multicolumn{5}{c}{$\rho=0.3$}\\
\multirow{3}{*}{100}
& Proposed & 59.8 & 56.3 & 55.7 & 53.1 & 53.1 & 99.7 & 99.8 &  99.7 & 100.0 & 99.9 \\
& CLX      & 13.9 &  8.9 &  8.1 &  6.9 &  6.6 & 51.5 & 45.7 &  38.3 &  31.9 & 27.2 \\
& LC       & 60.7 & 63.2 & 64.8 & 62.4 & 63.3 & 99.7 & 99.8 & 100.0 &  99.9 & 99.8 \\
\hline
\multirow{3}{*}{200}
& Proposed &  98.6 & 99.3 & 99.3 & 98.8 & 98.9 & 100.0 & 100.0 & 100.0 & 100.0 & 100.0 \\
& CLX      &  46.5 & 40.1 & 30.9 & 28.0 & 25.3 &  99.8 &  99.9 & 100.0 &  99.8 &  99.9 \\
& LC       &  98.6 & 99.3 & 99.0 & 99.1 & 98.9 & 100.0 & 100.0 & 100.0 & 100.0 & 100.0 \\
\hline
\end{tabular}

\vspace{1.5ex}
{\small
Note: This table reports the frequencies of rejection by each method under the null and alternative hypotheses based on $1000$ independent replications at the significance level $5\%$.}
\end{table}

For each covariance model, we generate samples independently from $N_p (\bzero, \bSigma^{*(i)})$ to evaluate the size, and use different covariance pairs described above to examine the power against dense and sparse alternatives. The empirical size and power are calculated  based on 1,000 replications at significance level $5\%$ and the results are reported in Tables \ref{tab: standardnormal}, \ref{tab: Model23} and \ref{tab: Model45}.

\begin{table}[H]
\small
\centering
\caption{Comparison of Empirical Size and Power (\%) for Models (ii) and (iii) }\label{tab: Model23}
\vspace{1ex}
\begin{tabular}{ccrrrrr|rrrrrrrrrrrrrrr}
  \hline
  && \multicolumn{5}{c}{Model (ii)} & \multicolumn{5}{c}{Model (iii)}\\
  \hline
 n&p& 100 & 200 & 500 & 800 & 1000 &  100 & 200 & 500 & 800 & 1000 \\
   \hline
 &&\multicolumn{10}{c}{Size}\\
\multirow{3}{*}{100}
& Proposed & 4.9 & 5.5 & 4.2 & 5.6 & 5.3 & 6.0 & 6.1 & 4.8 & 4.9 & 3.9 \\
& CLX      & 4.6 & 5.4 & 4.9 & 5.5 & 4.5 & 4.5 & 4.4 & 5.1 & 4.6 & 4.0\\
& LC       & 4.6 & 5.3 & 3.8 & 4.5 & 5.2 & 5.3 & 5.6 & 4.7 & 5.1 & 4.3\\
\hline
\multirow{3}{*}{200}
& Proposed & 6.5 & 5.4 & 4.1 & 3.8 & 4.3 & 6.3 & 6.5 & 4.8 & 4.1 & 4.9 \\
& CLX      & 4.5 & 4.3 & 5.8 & 4.0 & 4.3 & 4.3 & 6.5 & 4.1 & 3.8 & 4.8 \\
& LC       & 5.8 & 4.9 & 4.1 & 3.7 & 5.1 & 5.6 & 5.2 & 4.3 & 4.3 & 4.8 \\
\hline
 &&\multicolumn{10}{c}{Power under sparse alternative}\\
\multirow{3}{*}{100}
& Proposed & 98.4 & 96.1 & 87.5 & 85.3 & 79.8 & 98.1 & 95.7 & 88.1 & 82.3 & 81.3 \\
& CLX      & 98.8 & 97.7 & 92.3 & 90.2 & 85.9 & 98.7 & 97.5 & 91.3 & 88.0 & 86.6 \\
& LC       & 19.7 & 11.4 &  6.8 &  5.8 &  5.7 & 20.0 & 11.6 &  6.6 &  5.4 &  5.3 \\
\hline
\multirow{3}{*}{200}
& Proposed & 100.0 & 100.0 & 100.0 & 100.0 & 100.0 & 100.0 & 100.0 & 100.0 & 100.0 & 100.0\\
& CLX      & 100.0 & 100.0 & 100.0 & 100.0 & 100.0 & 100.0 & 100.0 & 100.0 & 100.0 & 100.0  \\
& LC       &  50.1 &  22.5 &   8.7 &   7.2 &   6.1 &  53.7 &  23.0 &  10.1 &   6.9 &   6.0 \\
\hline
 &&\multicolumn{10}{c}{Power under dense alternative}\\
\multirow{3}{*}{100}
& Proposed & 85.7 & 83.0 & 84.7 & 83.7 & 81.7 & 97.6 & 98.0 & 97.6 & 96.3 & 98.2 \\
& CLX      & 15.9 & 11.7 &  7.0 &  7.7 &  6.2 & 36.0 & 27.5 & 21.5 & 17.0 & 14.8 \\
& LC       & 88.5 & 87.7 & 89.6 & 89.2 & 89.8 & 97.9 & 98.5 & 98.5 & 97.4 & 99.1 \\
\hline
\multirow{3}{*}{200}
& Proposed & 100.0 & 100.0 & 100.0 & 100.0 & 100.0 & 100.0 & 100.0 & 100.0 & 100.0 & 100.0\\
& CLX      &  59.6 &  50.4 &  37.5 &  33.7 &  31.1 &  90.7 &  91.8 &  87.7 &  86.1 &  83.6\\
& LC       & 100.0 & 100.0 & 100.0 & 100.0 & 100.0 & 100.0 & 100.0 & 100.0 & 100.0 & 100.0\\
\hline
\end{tabular}

\vspace{1.5ex}
{\small
Note: This table reports the frequencies of rejection by each method under the null and alternative hypotheses based on $1000$ independent replications at the significance level $5\%$.}
\end{table}

\vspace{-2ex}
\begin{table}[H]
\small
\centering
\caption{Comparison of Empirical Size and Power (\%) Comparisons for Models (iv) and (v) }\label{tab: Model45}
\vspace{1ex}
\begin{tabular}{ccrrrrr|rrrrrrrrrrrrrrr}
  \hline
  && \multicolumn{5}{c}{Model (iv)} &  \multicolumn{5}{c}{Model (v)}\\
  \hline
 n&p& 100 & 200 & 500 & 800 & 1000 & 100 & 200 & 500 & 800 & 1000 \\
   \hline
 &&\multicolumn{10}{c}{Size}\\
\multirow{3}{*}{100}
& Proposed & 9.8 & 9.5 & 10.4 &  9.6 & 9.3 & 5.7 & 5.2 & 4.0 & 4.8 & 4.3 \\
& CLX      & 4.1 & 4.1 &  3.8 &  4.2 & 4.0 & 4.6 & 4.9 & 4.6 & 4.9 & 4.2  \\
& LC       & 9.5 & 9.3 & 10.7 & 10.3 & 9.3 & 5.4 & 5.2 & 4.7 & 4.6 & 3.7  \\
\hline
\multirow{3}{*}{200}
& Proposed & 10.1 & 10.8 & 9.0 & 10.1 & 8.2 & 6.3 & 6.0 & 3.6 & 4.3 & 4.4 \\
& CLX      &  3.2 &  4.5 & 3.0 &  3.4 & 4.8 & 5.1 & 4.0 & 3.7 & 4.6 & 4.3 \\
& LC       &  8.8 & 10.6 & 9.0 & 10.7 & 8.2 & 5.7 & 5.2 & 4.1 & 3.8 & 5.0 \\
\hline
 &&\multicolumn{10}{c}{Power under sparse alternative}\\
\multirow{3}{*}{100}
& Proposed & 97.6 & 96.4 & 88.1 & 84.8 & 81.5 &  99.9 & 85.2 & 78.9 & 72.5 & 86.7 \\
& CLX      & 98.8 & 98.1 & 92.4 & 89.3 & 86.5 & 100.0 & 90.0 & 83.5 & 77.8 & 90.9 \\
& LC       & 19.3 & 12.0 &  6.8 &  5.9 &  5.0 &  33.1 & 11.3 &  6.9 &  5.2 &  4.6 \\
\hline
\multirow{3}{*}{200}
& Proposed & 100.0 & 100.0 & 100.0 & 100.0 & 100.0 & 100.0 & 100.0 & 100.0 & 100.0 & 100.0\\
& CLX      & 100.0 & 100.0 & 100.0 & 100.0 & 100.0 & 100.0 & 100.0 & 100.0 & 100.0 & 100.0 \\
& LC       &  52.3 &  22.1 &   8.8 &   7.2 &   7.3 &  80.3 &  20.4 &   8.0 &   8.6 &   6.9\\
\hline
 &&\multicolumn{10}{c}{Power under dense alternative}\\
\multirow{3}{*}{100}
& Proposed & 84.1 & 89.7 & 92.2 & 95.5 & 96.8 & 100.0 & 100.0 & 100.0 & 100.0 & 100.0\\
& CLX      & 57.4 & 62.8 & 67.3 & 76.3 & 76.4 &  34.9 &  14.0 &   6.9 &   5.3 &   5.1\\
& LC       & 84.5 & 89.4 & 92.4 & 95.8 & 96.5 & 100.0 & 100.0 & 100.0 & 100.0 & 100.0 \\
\hline
\multirow{3}{*}{200}
& Proposed & 98.9 & 98.7 & 99.8 & 99.9 &  100.0 & 100.0 & 100.0 & 100.0 & 100.0 & 100.0\\
& CLX      & 88.6 & 90.3 & 95.6 & 97.1 &   98.0 &  94.2 &  52.0 &  12.8 &   8.8 &   6.9\\
& LC       & 99.1 & 98.9 & 99.9 & 99.8 &  100.0 & 100.0 & 100.0 & 100.0 & 100.0 & 100.0 \\
\hline
\end{tabular}

\vspace{1ex}
{\small
Note: This table reports the frequencies of rejection by each method under the null and alternative hypotheses based on $1000$ independent replications at the significance level $5\%$.}
\end{table}

The size and power comparisons from Tables \ref{tab: standardnormal}, \ref{tab: Model23} and \ref{tab: Model45} give us some intriguing findings:
\begin{itemize}
	\item[(1)] Under $H_0$, the sizes of all three tests are well retained close to the nominal level 0.05, except for Model (iv), in which both the LC test and our proposed test suffer from the size distortion, because of the violation of the test assumptions on covariance matrices.
	
	\vspace{-0.5ex}
	\item[(2)] As can be seen from Model (i), the CLX test is demonstrated to be powerful under the sparse alternative $H_s$, however, its performance is not satisfactory under the dense alternative. Even though in Models (ii)-(iv), the CLX test still has competitive powers, it fails with a decaying power as dimension grows in Model (v).
	
	\vspace{-0.5ex}
	\item[(3)] In the meantime, the LC test remains a high power under the dense alternative $H_d$, whereas performs poorly against the sparse alternative with a tendency of decaying as dimension $p$ grows large.
	
	\vspace{-0.5ex}
	\item[(4)] In comparison, our proposed Fisher's combined test exhibits competent results. Our proposed test performs as good as the CLX test under the sparse alternative, together with the comparable performance to the LC test when against the dense alternative.
\end{itemize}

In a summary, based on the simulation results in this section, we are able to say that the proposed Fisher test boost the power tremendously against more general alternatives, in the meanwhile, retaining the desired nominal significance level.

\section{Application to Gene-Set Testing}
We further demonstrate the power of our proposed test by applying the test to identify those sets of genes which potentially have significant differences in covariance matrices across different types of tumors. In biology, each gene does not work individually, but rather tends to function as groups to achieve complex biological tasks. Sets of genes are interpreted by Gene Ontology (GO) terms making use of the Gene Ontology system, in which genes are assigned to a set of predefined bins depending on their functional characteristics. The Gene Ontology covers three domains: biological process (BP), cellular component (CC) and molecular function (MF).

We consider the Acute Lymphoblastic Leukemia(ALL) data from the Ritz Laboratory at the Dana-Farber Cancer Institute (DFCI). The latest data is accessible at the ALL package (version 1.24.0) on \href{https://www.bioconductor.org/}{\color{blue}{Bioconductor}} website, including the original version published by \cite{chiaretti2004gene}. The ALL dataset consists of microarrays expression measures of 12,625 probes on Affymetrix chip series HG-U95Av2 for 128 different individuals with acute lymphoblastic leukemia, which is a type of blood cancer in that bone marrow affects white blood cells.
Based on the type of lymphocyte that the leukemia cells come from, the disease is classified into subgroups of T-cell ALL and B-cell ALL. In our study, we focus on a subset of the ALL data of 79 patients with the B-cell ALL. We are interested in two types of B-cell tumors: BCR/ABL and NEG, with sample sizes being 37 and 42 respectively.

Let us consider $K$ gene sets $S_1,\cdots, S_K$, and $\bSigma_{1S_k}$ and $\bSigma_{2S_k}$ be the covariance matrices of two types of tumors respectively. The null hypotheses we are interested are

$$H_{0,category}: \bSigma_{1S_k} = \bSigma_{2S_k},\quad k=1,\cdots,K$$
where $category \in \{BP, CC, MF\}$ because we classify the gene sets into three different GO categories and shall test each GO category separately.

\begin{figure}[H]
\centering
\caption{Histograms of the Dimension of Gene-sets for Three GO Categories}\label{fig: plotgeneset}
\vspace{1ex}
\includegraphics[width=\textwidth, height = 2.3in]{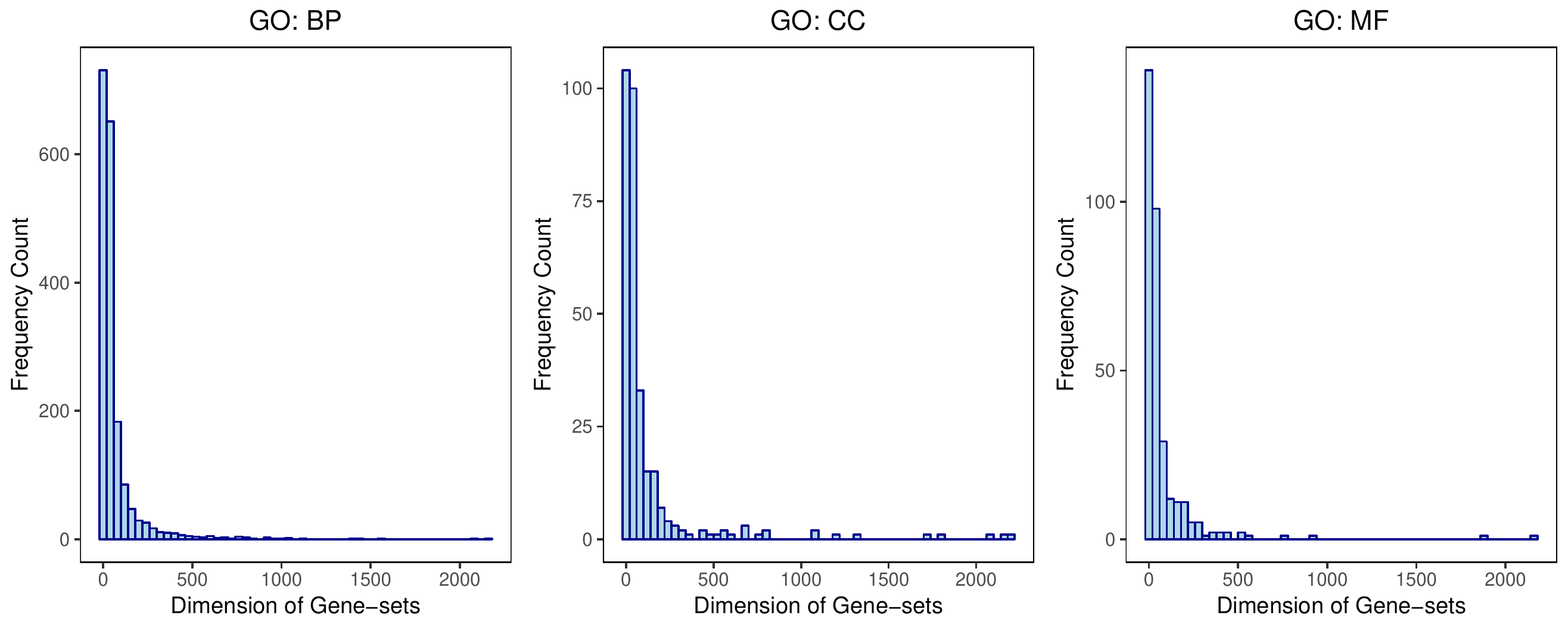}
\end{figure}

To control the computational costs, we first perform a pre-screening procedure following the same criteria as in \cite{dudoit2008multiple} by choosing those probes that satisfy (i) the fluorescence intensities greater than 100 (absolute scale) for at least 25\% of the 79 cell samples; (ii) the interquartile range (IQR) of the fluorescence intensities for the 79 cell samples greater than 0.5 (log base 2 scale). The preliminary gene-filtering retains 2,391 probes. After that we then identify those GO terms annotating at least 10 of the 2,391 filtered probes, which gives us 1849 unique GO terms in BP category, 306 in CC and 324 in MF for further analysis. Table \ref{tab: summary-gene-set} and Figure \ref{fig: plotgeneset} summarize the dimension of gene-sets contained in each category.

\vspace{-2ex}
\begin{table}[H]
\small
\centering
\caption{Summary of the Dimension of Gene-sets for Three GO Categories}\label{tab: summary-gene-set}
\vspace{1ex}
\begin{tabular}{c|c|ccccccc}
\hline
GO Category & Total number & Min & 1st-Quantile & Median & 3rd-Quantile & Max  \\ 
 \hline
BP & 1849 & 10 & 15 & 27 & 62 & 2153  \\ \hline 
CC &  306 & 10 & 17 & 32  & 85 & 2181 \\ \hline 
MF &  324 & 10 & 14 & 26 & 68 & 2148  \\ 
\hline
\end{tabular}
\end{table}

\vspace{-3ex}

\begin{figure}[H]
\centering
\caption{Boxplots of the LC and CLX Test Statistics for Three GO Categories}\label{fig: plotstat}
\vspace{1ex}
\includegraphics[width=5in, height=3in]{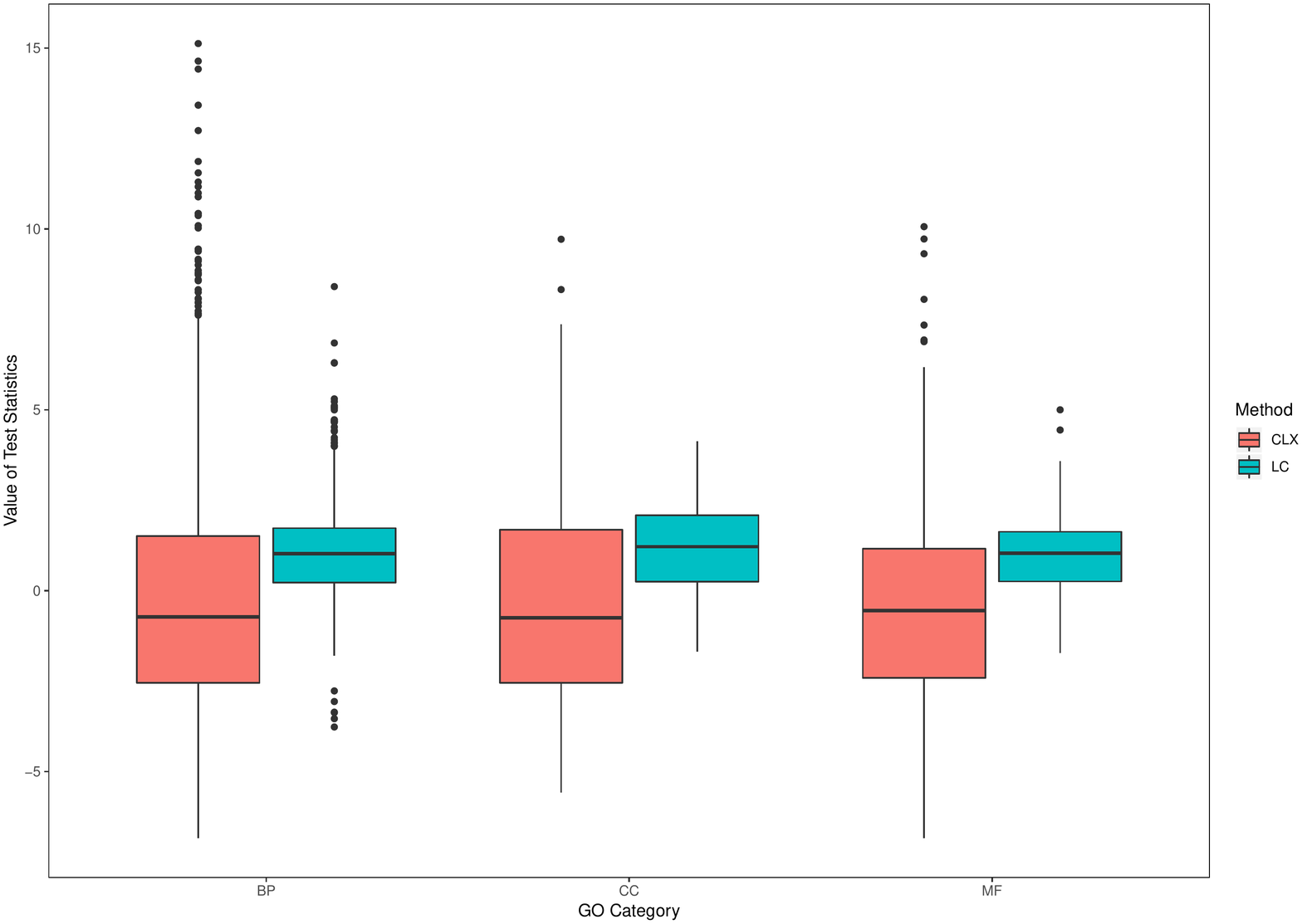}
\end{figure}

We first take a look at the performance of the CLX test and the LC test. Figure \ref{fig: plotstat} displays boxplots of both test statistics. It can be observed that test statistics have quite different magnitudes, indicating difficulty in the approach of weighted summation combination of the two statistics.

\vspace{-2ex}

We then apply our proposed Fisher's method to test the hypothesis, together with comparisons to the CLX and LC tests. We also compare our test with the natural Bonferroni combination. The test outcomes are reported in Table \ref{tab: realdata}, with nominal level $\alpha=0.05$ for each test. Furthermore, in order to control the false discovery rate (FDR), we apply the Benjamini-Hochberg (BH) procedure \citep{benjamini1995controlling} to each GO category, and the results are listed in Table \ref{tab: realdataBH}, with nominal level $\alpha=0.05$ for every category.

\begin{table}[H]
\small
\centering
\caption{Gene-Set Testing Results at the Nominal Level $\alpha=0.05$} \label{tab: realdata}
\vspace{1ex}
\begin{tabular}{c|c|ccccc}
\hline
\multirow{2}{*}{GO Category} & Total number & \multicolumn{4}{c}{Number of Significant Gene-sets} \\
\cline{3-6}
& of Gene-sets & \hspace{0.3cm} CLX \hspace{0.3cm} & \hspace{0.3cm} LC \hspace{0.4cm} & Bonferroni & Proposed \\
 \hline
BP & 1849 & 297 & 505& 451 & 615 \\ \hline
CC &   306 &  52 & 111 & 96 & 116 \\ \hline
MF &   324 &  38 & 78 & 61 & 96\\
\hline
\end{tabular}
\end{table}

\begin{table}[H]
\small
\centering
\caption{Gene-Set Testing Results with the FDR Control at $\alpha=0.05$}\label{tab: realdataBH}
\vspace{1ex}
\begin{tabular}{c|c|ccccc}
\hline
\multirow{2}{*}{GO Category} & Total number & \multicolumn{4}{c}{Number of Significant Gene-sets} \\
\cline{3-6}
& of Gene-sets & \hspace{0.3cm} CLX \hspace{0.3cm} & \hspace{0.3cm} LC \hspace{0.4cm} & Bonferroni & Proposed \\
 \hline
BP &1849  & 0 & 126 & 81 & 254 \\ \hline
CC & 306 & 0 & 55 & 24 &68 \\ \hline
MF & 324 & 0 & 20 & 4 & 26\\
\hline
\end{tabular}
\end{table}

As shown in Table \ref{tab: realdataBH}, our proposed test identifies much more significant gene-sets than the other methods. The LC identifies a few while the Bonferroni test identifies fewer significant gene-sets than the LC test does. This illustrates that the Bonferroni test is relatively conservative, which is consistent with what we expect. Unfortunately, the CLX test fails to declare any significance after we control the FDR using BH procedure. This is possibly because the signals in the differences are not strong enough for the CLX test to detect.

Biological evidence supports that such improvement is quite meaningful and very helpful in cancer research. To clarify this, we further investigate those gene-sets that are not declared significant by the CLX and LC tests but are identified by our proposed Fisher test. Taking the GO term ``GO:0005905" as an example, it refers to the clathrin-coated pit which functions in the cellular component (CC) gene ontology category. Protein evidence by \cite{ezkurdia2014multiple} confirms that the clathrin-coated pit works with several protein-coding genes, such as CLTCL1, PICALM, etc., that are closely related to human cancers. We also take a deep look at ``GO:0035259", the glucocorticoid receptor binding, in the molecular function (MF) gene ontology category. Many genes contribute to this gene-set, among them, we pay special attention to STAT3, a protein-coding gene which plays an important role in the immune system by transmitting signals for the maturation of immune system cells, especially T-cells and B-cells. Researchers have observed that STAT3 gene mutations are highly correlated with cancers, especially blood cancers \citep{hodge2005role,jerez2012stat3,haapaniemi2015autoimmunity,milner2015early}. In a short summary, our proposed test incorporates the information from the CLX statistic, which successfully enhances the power over the LC test, even though the LC test itself may not declare any significance.

\section{Conclusion}
This paper studies the fundamental problem of testing high-dimensional covariance matrices. Unlike the existing quadratic form statistics, maximum form statistics, and their weighted combination, we provide a new perspective to exploit the full potential of quadratic form statistics and maximum form statistics. We propose a scale-invariant and computationally efficient power enhancement test based on Fisher's method to combine their respective $p$-values. Theoretically, after deriving their joint limiting null distribution, we prove that the proposed combination method retains the correct asymptotic size and boosts the power against more general alternatives. Numerically, we demonstrate the finite-sample properties in simulation studies and the practical relevance through an empirical study on gene-set testing problem.

It is still an open question to relax the Gaussian assumption when deriving the asymptotic joint distribution of quadratic form statistics and maximum form statistics in the two-sample covariance tests. There are several potential directions to relax the Gaussian assumption. For instance, we may use the semiparametric Gaussian copula distribution \citep{liu2012high,xue2012regularized} and study the nonparametric tests. Alternatively, we may use the Gaussian approximation theory to bridge this gap. We will leave this open question for future work.

{
\bibliographystyle{agsm}
\bibliography{cov2testref}
}

\end{document}